\documentclass[11pt,leqno]{amsart}
\usepackage{amsmath, amsthm, amscd, amsfonts, amssymb, graphicx, color}
\usepackage[bookmarksnumbered, plainpages]{hyperref}

\newtheorem{thm}{Theorem}[section]

\newtheorem{prop}[thm]{Proposition}

\numberwithin{equation}{section}

\usepackage{latexsym}
\begin{document}

\title{Zeta and normal zeta functions for a subclass of  space groups}


\author{Hermina ALAJBEGOVI\'C \and  Muharem AVDISPAHI\'C }



\maketitle

\begin{abstract}
We calculate zeta and normal zeta functions of space groups with the
point group isomorphic to the cyclic group of order 2. The obtained results are applied to determine the number of subgroups, resp. normal subgroups, of a
given index for each of these groups.\\\\
\\
{\bf Mathematics Subject Classification (2010) 11M41, 20H15}

\end{abstract}

\section{\bf Introduction}

The concepts of the zeta and normal zeta function of a group were introduced by Smith, Segal and
Grunewald in \cite {GSS} and \cite {S}. The zeta function of a group
$G$ is defined as ${\zeta _G}(s) = \sum\limits_{n \in \mathbb {N}}{{a_n}(G){n^{ - s}}} $, where ${a_n}(G)$ denotes the number of
subgroups of index $n$ in $G$. The normal zeta
function of a group $G$ is given by $\zeta _G^ \triangleleft (s) =
\sum\limits_{n \in \mathbb {N}} {{c_n}(G){n^{ - s}}} $, where ${c_n}(G)$ is
the number of normal subgroups of index $n$ in $G$. These functions provide a useful tool for studying the relationship
 between the asymptotic behavior of the sequences ${a_n(G)}$, resp. ${c_n(G)}$, and the structure of $G$.  

\mbox{}If a group $G$ has a polynomial subgroup growth, i.e.  if ${{a_n}(G)\leq n^{k}}$ for some $k$  and for all  $n\in  N$,
then ${\zeta _G}(s)$ has a non - trivial domain of convergence. In particular, if
$G$ is residually finite nilpotent group then ${a_n}(G)$ grows
polynomially and ${\zeta _G}(s) = \sum\limits_{n \in \mathbb {N}}
{{a_n}(G){n^{ - s}}} $ satisfies an Euler product formula ${\zeta _G}(s) =
\prod\limits_{p \in P} {{\zeta _{G,p}}(s)} $, where ${\zeta
_{G,p}}(s)$ counts only subgroups of $p$- power index and $P$
denotes the set of all primes.  
 
M. P. F. du Sautoy, J.
J. McDermott and G. C. Smith \cite {duSMS} proved the following
theorem.

\begin{thm} Let $G$ be a finite extension of a free abelian group of
finite rank. Then ${\zeta _G}(s)$ and $\zeta _G^ \triangleleft (s)$
can be extended to meromorphic functions on the whole complex plane.
\end{thm}

Lubotzky and du Sautoy \cite {duSL} established a
functional equation ${\left. {{{ \zeta }_{G,p}}(s)} \right|_{p
\to {p^{ - 1}}}} = {( - 1)^n}{p^{as + b}}{ \zeta _{G,p}}(s)$
satisfied by the local factors of the zeta function of a group for some torsion-free nilpotent groups and 
appropriate $a,b,n\in \mathbb {N}$. 
Here $p \to {p^{ - 1}}$ denotes a formal inversion of the local parameter $p$. In this regard, see also \cite {C}, \cite {CB}.

\mbox{}There are relatively few explicit expressions known for
zeta functions of groups. John J. McDermott calculated the zeta functions of the seventeen
plane crystallographic groups in \cite {M}.

\mbox{}A space group 
represents a description of the symmetry of a crystal. A crystallographic group $G$ contains a translation subgroup $T$ which
consists of all elements of the group corresponding to translations
of the pattern involved. The factor group ${G
\mathord{\left/{\vphantom {G T}}\right.\kern-\nulldelimiterspace}
T}$ is known as the point group of $G$ and is denoted by $P$. Group
$G$ is a finite extension of $T$ by $P$, since $P$ is finite. 
 
\mbox{}In this paper, we calculate zeta and normal zeta functions of space groups with the 
point group isomorphic to the cyclic group of order 2. There are eight such groups: ${P{\bar 1}}$,
  ${P2}$, ${P{{2_1}}}$,  ${C2}$, ${Pm}$,  ${Pc}$, ${Cm}$ and  ${Cc}$ \cite {L}.
After stating results in Section 2, we provide the application of these results to compute the number of subgroups 
of a given index for each group in Section 3 and Section 4. We describe the method of
 proof in Section 5  and provide full details in case of the group C2, as a sample.

\section{\bf Results}

For a sake of bravity, the following notation for translates of the
Riemann zeta
function is used in the sequel:\\ ${\zeta_k(s)=\zeta(s-k)}$, i.e., ${\zeta_2(s)=\zeta(s-2)}$.\\

\begin{thm} Zeta functions of space groups with the point group isomorphic to the cyclic group of order 2 read as follows 

$\zeta_{P{\bar 1}}(s)= {\zeta_1(s)\zeta_2(s)\zeta_3(s) + {2^{ - s}}\zeta(s)\zeta_1(s)\zeta_2(s)}$ 

$\zeta_ {P2}(s)= {(1 + {2^{ - s + 3}})\zeta (s)\zeta_1(s)\zeta_2(s)} $

$ \zeta_{P{{2_1}}}(s)={\zeta (s)\zeta_1(s)\zeta_2(s)} $  

$\zeta_{C2}(s) = {(1 + {2^{ - 2s + 3}})\zeta (s)\zeta_1(s)\zeta_2(s)}$ 

$\zeta_ {Pm}(s)={(1 + 9 \cdot {2^{ - s}} + 6 \cdot {2^{ - 2s}})\zeta
(s)\zeta_1(s)\zeta_1(s) + {2^{ - s}}\zeta (s)\zeta_1(s)\zeta_2(s)} $

$\zeta_{Pc}(s)= {(1 + {2^{ - s}} - 2 \cdot {2^{ - 2s}})\zeta (s)\zeta_1(s)\zeta_1(s) + {2^{ - s}}\zeta (s)\zeta_1(s)\zeta_2(s)}$ 

$\zeta_{Cm}= {(1 + {2^{ - s}}+6 \cdot {2^{ - 2s}}+8 \cdot {2^{ -
3s}})\zeta (s)\zeta_1(s)\zeta_1(s) + {2^{ - s}}\zeta
(s)\zeta_1(s)\zeta_2(s)}  $ 

$\zeta_{Cc}(s)={(1 - 3 \cdot{2^{ - s}}+10 \cdot {2^{ - 2s}}-8 \cdot {2^{
- 3s}})\zeta (s)\zeta_1(s)\zeta_1(s) + {2^{ - s}}\zeta(s)\zeta_1(s)\zeta_2(s)} $

\end{thm}

\begin{thm} Normal zeta functions of space groups with the point
group isomorphic to the cyclic group of order 2 are given by

$\zeta_{P{\bar 1}}^\triangleleft (s)={1 + 14 \cdot {2^{ - s}} + 28 \cdot
{2^{ - 2s}}+8 \cdot {2^{ - 3s}} + {2^{-s}}\zeta(s)\zeta_1(s)\zeta_2(s)} $

$\zeta_{P2}^
\triangleleft (s)={(1 + 13 \cdot {2^{ - s}} + 22 \cdot {2^{ - 2s}} + 4 \cdot {2^{ -
3s}})\zeta (s) + (3 \cdot {2^{ -2s}} + {2^{ - s}})\zeta (s)\zeta (s)\zeta_1(s)}$

$\zeta_{P{{2_1}}}^
\triangleleft (s)={(1 +5\cdot{2^{ - s}} - 2 \cdot {2^{ - 2s}} - 4 \cdot {2^{ - 3s}})\zeta
(s) + ({2^{ - s}} + 3 \cdot {2^{ -  2s}})\zeta (s)\zeta (s)\zeta_1(s)}$

$\zeta_{C2}^
\triangleleft (s)={(2 \cdot {2^{ - 2s}} + 5 \cdot {2^{ - s}} + 1)\zeta (s) + {2^{ -
s}} \cdot (1 - {2^{ - s}} + 4 \cdot {2^{ -
2s}})\zeta (s)\zeta (s)\zeta_1(s)}$

$\zeta_{Pm}^\triangleleft (s) ={(1 + 11 \cdot {2^{ - s}} + 12 \cdot {2^{ - 2s}})\zeta
(s)\zeta_1(s) + {2^{ - s}}(1 + 3 \cdot {2^{ -
s}})\zeta (s)\zeta (s)\zeta_1(s)}$

$\zeta_{Pc}^\triangleleft (s)={(1 + 3 \cdot {2^{ - s}} - 4 \cdot {2^{ - 2s}})\zeta
(s)\zeta_1(s) + (3 \cdot {2^{ - 2s}} + {2^{ -
s}})\zeta (s)\zeta (s)\zeta_1(s)} $

$\zeta_{Cm}^
\triangleleft (s)= {(3 \cdot {2^{ - s}} + 1)\zeta
(s)\zeta_1(s) + {2^{ - s}} \cdot (1 - {2^{ - s}} + 4 \cdot {2^{ -
2s}})\zeta(s)\zeta (s)\zeta_1(s)}  $

 $\zeta _{Cc}^
\triangleleft (s)= {(1 - {2^{ - s}})\zeta (s)\zeta_1(s) +
(4 \cdot {2^{ -
3s}} - {2^{ - 2s}} + {2^{ - s}})\zeta (s)\zeta (s)\zeta_1(s)}$
\end{thm}

\section{\bf Applications: Subgroups of a given index}

In this and the following section,  $d(n)$ denotes the number of all positive divisors of a positive integer $n$ and $\sigma (n)$ denotes the sum of all positive divisors for a positive integer $n$, i. e. $\sigma (n) = \sum\limits_{\left. l \right|n} l .$

\begin{prop} The number of all  subgroups of  index $n$ in the group $P\bar 1$ is given by the following expressions
\begin{enumerate}
  \item  if $n$ is even, \[{a_n} = n\sum\limits_{\left. l \right|n} {l \cdot \sigma (l) + } \sum\limits_{\left. l \right|\left( {\frac{n}{2}} \right)} {l \cdot \sigma (l)} \] 
\item	if $n$ is odd, ${a_n} = n\sum\limits_{\left. l \right|n} {l \cdot \sigma (l)},$

\end{enumerate}
In particular, ${a_p} = {p^3} + {p^2} + p$ for every odd prime $p$. 
\end{prop}

\begin {prop} The number of all  subgroups of  index $n$ in group  $P2$ reads:
\begin{enumerate}
  \item  if $n$ is even, \[{a_n}  = \sum\limits_{\left. l \right|n} {l \cdot \sigma (l)}  + 8 \cdot \sum\limits_{\left. l \right|\left( {\frac{n}{2}} \right)} {l \cdot \sigma (l)} \] 
\item  if $n$ is odd, ${a_n} = \sum\limits_{\left. l \right|n} {l \cdot \sigma (l)},$
\end{enumerate}
In particular, if $p$ is an odd prime, then ${a_p} = {p^2} + p + 1.$ 
\end {prop}

\begin {prop} The number of all  subgroups of  index $n$ in group  ${P{{2_1}}}$ is 
${a_n} =  \sum\limits_{\left. l \right|n} {l \cdot \sigma (l)}.$  In particular, ${a_p} = {p^2} + p + 1$ for every odd prime $p$.
\end{prop}

\begin {prop}   The number of all  subgroups of  index $n$ in the  group  $C2$  is

\[{a_n} = \left\{ \begin{gathered}
  \sum\limits_{\left. l \right|n} {l \cdot \sigma (l)} ,\,\,\,(n \equiv 1 \vee n \equiv 2 \vee n \equiv 3)\,(\bmod 4) \hfill \\
  \sum\limits_{\left. l \right|n} {l \cdot \sigma (l)}  + 8 \cdot \sum\limits_{\left. l \right|\left( {\frac{n}{4}} \right)} {l \cdot \sigma (l)},  n \equiv 0\,(\bmod 4) \hfill \\ 
\end{gathered}  \right.\]
In particular, if $p$ is an odd prime, then ${a_p} = 1 + p + {p^2}.$ 
\end{prop}

\begin {prop}  The number of all  subgroups of  index $n$ in the group  $Pm$ is as follows:
\begin{enumerate}

  \item if $n$  is even,\\
${a_n} = \left\{ \begin{gathered}
  \sum\limits_{\left. l \right|n} {l \cdot d(l)}  + 9 \cdot \sum\limits_{\left. l \right|\left( {\frac{n}{2}} \right)} {l \cdot d(l)}  + \sum\limits_{\left. l \right|\left( {\frac{n}{2}} \right)} {l \cdot \sigma (l)} ,n \equiv 2\,(\bmod 4) \hfill \\
  \sum\limits_{\left. l \right|n} {l \cdot d(l)}  + 9  \sum\limits_{\left. l \right|\left( {\frac{n}{2}} \right)} {l \cdot d(l)}  + 6  \sum\limits_{\left. l \right|\left( {\frac{n}{4}} \right)} {l \cdot d(l)}  + \sum\limits_{\left. l \right|\left( {\frac{n}{2}} \right)} {l \cdot \sigma (l)} ,n \equiv 0\,(\bmod 4) \hfill \\ 
\end{gathered}  \right.$
\item if $n$ is odd, ${a_n}= \sum\limits_{\left. l \right|n} {l \cdot d(l)}.$ 

\end{enumerate}
   
In particular, if $p$ is an odd prime, then ${a_p} = 2p + 1.$
\end{prop}

\begin{prop}  The number of all  subgroups of  index $n$ in the group  $Pc$ is given by:
\begin{enumerate}

  \item  if $n$  is even, \\
${a_n} = \left\{ \begin{gathered}
  \sum\limits_{\left. l \right|n} {l \cdot d(l)}  + \sum\limits_{\left. l \right|\left( {\frac{n}{2}} \right)} {l \cdot d(l)}  + \sum\limits_{\left. l \right|\left( {\frac{n}{2}} \right)} {l \cdot \sigma (l)} ,n \equiv 2\,(\bmod 4) \hfill \\
  \sum\limits_{\left. l \right|n} {l \cdot d(l)}  + \sum\limits_{\left. l \right|\left( {\frac{n}{2}} \right)} {l \cdot d(l)}  - 2 \cdot \sum\limits_{\left. l \right|\left( {\frac{n}{4}} \right)} {l \cdot d(l)}  + \sum\limits_{\left. l \right|\left( {\frac{n}{2}} \right)} {l \cdot \sigma (l)} ,n \equiv 0\,(\bmod 4) \hfill \\ 
\end{gathered}  \right.$
\item	if $n$ is odd, ${a_n} = \sum\limits_{\left. l \right|n} {l \cdot d(l)}.$ 

\end {enumerate}
In particular, if $p$ is an odd prime, then ${a_p} = 2p + 1.$
\end{prop}

\begin{prop}  The number of all  subgroups of  index $n$ in the group  $Cm$ is:
\begin{enumerate}	
\item if $n$ is even, \\
$a_n = \left\{ \begin{gathered}
  \sum\limits_{\left. l \right|n} {l \cdot d(l)}  + \sum\limits_{\left. l \right|\left( {\frac{n}{2}} \right)} {l \cdot d(l)}  + \sum\limits_{\left. l \right|\left( {\frac{n}{2}} \right)} {l \cdot \sigma (l)} ,(n \equiv 2\, \vee n \equiv 6\,)(\bmod 8) \hfill \\
  \sum\limits_{\left. l \right|n} {l \cdot d(l)}  + \sum\limits_{\left. l \right|\left( {\frac{n}{2}} \right)} {l \cdot d(l)}  + 6 \cdot \sum\limits_{\left. l \right|\left( {\frac{n}{4}} \right)} {l \cdot d(l)}  + \sum\limits_{\left. l \right|\left( {\frac{n}{2}} \right)} {l \cdot \sigma (l)} ,n \equiv 4\,(\bmod 8) \hfill \\
  \sum\limits_{\left. l \right|n} {l  d(l)}  + \sum\limits_{\left. l \right|\left( {\frac{n}{2}} \right)} {l d(l)}+ 6\sum\limits_{\left. l \right|\left( {\frac{n}{4}} \right)} {l d(l)}+8\sum\limits_{\left. l \right|\left( {\frac{n}{8}} \right)} {l d(l)}+\sum\limits_{\left. l \right|\left( {\frac{n}{2}} \right)} {l \sigma (l)} ,n \equiv 0\,(\bmod 8)\, \hfill \\ 
\end{gathered}  \right.$
\item	if $n$ is odd, ${a_n} = \sum\limits_{\left. l \right|n} {l \cdot d(l)}.$ 

\end {enumerate}
In particular, if $p$ is an odd prime, then ${a_p} = 2p + 1.$
\end{prop}

\begin{prop}  The number of all  subgroups of  index $n$ in the group  $Cc$ is the following:
 \begin{enumerate}
\item 	if $n$ is even,\\
${a_n} = \left\{ \begin{gathered}
  \sum\limits_{\left. l \right|n} {l \cdot d(l)}  - 3 \cdot \sum\limits_{\left. l \right|\left( {\frac{n}{2}} \right)} {l \cdot d(l)}  + \sum\limits_{\left. l \right|\left( {\frac{n}{2}} \right)} {l \cdot \sigma (l)} ,(n \equiv 2\, \vee n \equiv 6)\,(\bmod 8) \hfill \\
  \sum\limits_{\left. l \right|n} {l \cdot d(l)}  - 3 \cdot \sum\limits_{\left. l \right|\left( {\frac{n}{2}} \right)} {l \cdot d(l)}  + 10 \cdot \sum\limits_{\left. l \right|\left( {\frac{n}{4}} \right)} {l \cdot d(l)}  + \sum\limits_{\left. l \right|\left( {\frac{n}{2}} \right)} {l \cdot \sigma (l)} ,n \equiv 4\,(\bmod 8) \hfill \\
  \sum\limits_{\left. l \right|n} {l d(l)}  - 3 \sum\limits_{\left. l \right|\left( {\frac{n}{2}} \right)} {l d(l)}  + 10 \sum\limits_{\left. l \right|\left( {\frac{n}{4}} \right)} {l d(l)}  - 8 \sum\limits_{\left. l \right|\left( {\frac{n}{8}} \right)} {l d(l)}  + \sum\limits_{\left. l \right|\left( {\frac{n}{2}} \right)} {l \sigma (l)} ,n \equiv 0\,(\bmod 8) \hfill \\ 
\end{gathered}  \right.$

\item	if $n$ is odd, ${a_n} = \sum\limits_{\left. l \right|n} {l \cdot d(l)}.$
\end {enumerate} 
In particular, if $p$ is an odd prime, then ${a_p} = 2p + 1.$
\end{prop}

\section{\bf Applications: Normal subgroups of a given index}

\begin{prop} The number of all normal subgroups of index $n$ in the group $P\bar 1$ reads:
\begin{enumerate}
\item  $c_1 = 1$,
\item if $n$ is odd and $ n\neq 1$, then $c_n = 0$,
\item if $n$ is even, $c_n = \left\{ {\begin{array}{*{20}{c}}
   {15,\,\,\,\,\,\,\,\,\,\,\,\,\,\,\,\,\,\,\,\,\,\,\,\,\,\,\,\,\,\,\,\,\,\,\,\,\,\,\,\,\,\,\,\,\,\,\,\,\,\,\,\,\,\,\,\,\,\,\,\,\,\,n = 2,} \\ 
  {35,\,\,\,\,\,\,\,\,\,\,\,\,\,\,\,\,\,\,\,\,\,\,\,\,\,\,\,\,\,\,\,\,\,\,\,\,\,\,\,\,\,\,\,\,\,\,\,\,\,\,\,\,\,\,\,\,\,\,\,\,\,\,n = 4,} \\ 
  {43,\,\,\,\,\,\,\,\,\,\,\,\,\,\,\,\,\,\,\,\,\,\,\,\,\,\,\,\,\,\,\,\,\,\,\,\,\,\,\,\,\,\,\,\,\,\,\,\,\,\,\,\,\,\,\,\,\,\,\,\,\,\,n = 8,} \\ 
  {\sum\limits_{\left. l \right|\left( {\frac{n}{2}} \right)} {l \cdot \sigma (l)} ,\,\,\,\,\,n \equiv 0(\bmod 2) \wedge n \ne 1,2,4,8\,} 
\end{array}} \right.$
\end{enumerate}
\end{prop}

\begin{prop}  The number of all normal subgroups of index $n$ in the group $P2$ is:
\begin {enumerate}	
\item	if $n$  is even, 
$c_n = \left\{ \begin{gathered}
  40 + 3 \cdot \sum\limits_{l\left| {\left( {\frac{n}{4}} \right)} \right.} {\sigma (l)}  + \sum\limits_{l\left| {\left( {\frac{n}{2}} \right)} \right.} {\sigma (l)} ,\,\,\,n \equiv 0(\bmod 8), \hfill \\
  14 + \sum\limits_{l\left| {\left( {\frac{n}{2}} \right)} \right.} {\sigma (l)} ,\,\,(n \equiv 2 \vee \,n \equiv 6)(\bmod 8) \hfill \\
  36 + 3 \cdot \sum\limits_{l\left| {\left( {\frac{n}{4}} \right)} \right.} {\sigma (l)}  + \sum\limits_{l\left| {\left( {\frac{n}{2}} \right)} \right.} {\sigma (l)} ,\,\,\,\,n \equiv 4(\bmod 8) \hfill \\ 
\end{gathered}  \right.$
\item  if $n$ is odd, $c_n = 1.$
\end{enumerate}
\end {prop}

\begin{prop}  The number of all normal subgroups of index $n$ in the group $P{2_1}$ is given by:
\begin {enumerate}
\item	if $n$  is even,
$c_n = \left\{ \begin{gathered}
  3 \cdot \sum\limits_{l\left| {\left( {\frac{n}{4}} \right)} \right.} {\sigma (l)}  + \sum\limits_{l\left| {\left( {\frac{n}{2}} \right)} \right.} {\sigma (l)} ,\,\,\,\,n \equiv 0(\bmod 8),\,\, \hfill \\
  6 + \sum\limits_{l\left| {\left( {\frac{n}{2}} \right)} \right.} {\sigma (l)} ,\,\,\,(n \equiv 2 \vee \,n \equiv 6)(\bmod 8), \hfill \\
  4 + 3 \cdot \sum\limits_{l\left| {\left( {\frac{n}{4}} \right)} \right.} {\sigma (l)}  + \sum\limits_{l\left| {\left( {\frac{n}{2}} \right)} \right.} {\sigma (l)} ,\,\,n \equiv 4(\bmod 8) \hfill \\ 
\end{gathered}  \right.$
\item  	if $n$ is odd, $c_n = 1.$ 
\end{enumerate}
\end{prop}

\begin {prop}  The number of all normal subgroups of index $n$ in the group $C2$ is:
\begin{enumerate}
\item if $n$ is even,
  \[c_n = \left\{ \begin{gathered}
  6 + \sum\limits_{l\left| {\left( {\frac{n}{2}} \right)} \right.} {\sigma (l)} ,\, (n \equiv 2\, \vee n \equiv 6\,)(\bmod 8) \hfill \\
  8 - \sum\limits_{l\left| {\left( {\frac{n}{4}} \right)} \right.} {\sigma (l)}  + \sum\limits_{l\left| {\left( {\frac{n}{2}} \right)} \right.} {\sigma (l)} ,n \equiv 4\,(\bmod 8) \hfill \\
  8 + 4 \cdot \sum\limits_{l\left| {\left( {\frac{n}{8}} \right)} \right.} {\sigma (l)}  - \sum\limits_{l\left| {\left( {\frac{n}{4}} \right)} \right.} {\sigma (l)}  + \sum\limits_{l\left| {\left( {\frac{n}{2}} \right)} \right.} {\sigma (l)} ,n \equiv 0\,(\bmod 8)\, \hfill \\ 
\end{gathered}  \right.\] 
  \item if $n$ is odd, $c_n  = 1.$ 
\end {enumerate}
\end{prop}

\begin {prop}  The number of all normal subgroups of index $n$ in the group $Pm$ is:
	\begin{enumerate}
  \item  if $n$ is even, \[c_n = \left\{ \begin{gathered}
  \sigma (n) + 11 \cdot \sigma \left( {\frac{n}{2}} \right) + \sum\limits_{l\left| {\left( {\frac{n}{2}} \right)} \right.} {\sigma (l)} ,n \equiv 2\,(\bmod 4) \hfill \\
  \sigma (n) + 11 \cdot \sigma \left( {\frac{n}{2}} \right) + 12 \cdot \sigma \left( {\frac{n}{4}} \right) + \sum\limits_{l\left| {\left( {\frac{n}{2}} \right)} \right.} {\sigma (l)}  + 3 \cdot \sum\limits_{l\left| {\left( {\frac{n}{4}} \right)} \right.} {\sigma (l),\,} n \equiv 0\,(\bmod 4) \hfill \\ 
\end{gathered}  \right.\] 
\item if $n$ is odd, ${c_n}= \sigma (n).$
\end{enumerate}
In particular, if $p$ is an odd prime, then  ${c_p = p + 1}.$ 
\end{prop}

\begin{prop} The number of all normal subgroups of index $n$ in the group $Pc$ reads:
\begin{enumerate}
  \item  if $n$ is even, 
\[c_n  = \left\{ \begin{gathered}
  \sigma (n) + 3 \cdot \sigma \left( {\frac{n}{2}} \right) + \sum\limits_{l\left| {\left( {\frac{n}{2}} \right)} \right.} {\sigma (l)} ,\,\,\,n \equiv 2\,(\bmod 4) \hfill \\
  \sigma (n) + 3 \cdot \sigma \left( {\frac{n}{2}} \right) - 4 \cdot \sigma \left( {\frac{n}{4}} \right) + \sum\limits_{l\left| {\left( {\frac{n}{2}} \right)} \right.} {\sigma (l)}  + 3 \cdot \sum\limits_{l\left| {\left( {\frac{n}{4}} \right)} \right.} {\sigma (l)} ,\,n \equiv 0\,(\bmod 4) \hfill \\ 
\end{gathered}  \right.\] 
\item  if $n$ is odd, $c_n = \sigma (n).$

\end{enumerate} 
In particular, if $p$ is an odd prime, then ${c_p  = p + 1}.$ 
\end{prop}

\begin{prop} The number of all normal subgroups of index $n$ in the group $Cm$ is given by:
     \begin{enumerate}
\item if $n$ is even, \\
$c_n  = \left\{ \begin{gathered}
  \sigma (n) + 3 \cdot \sigma \left( {\frac{n}{2}} \right) + \sum\limits_{l\left| {\left( {\frac{n}{2}} \right)} \right.} {\sigma (l)} ,(n \equiv 2 \vee n \equiv 6)\,(\bmod 8) \hfill \\
  \sigma (n) + 3 \cdot \sigma \left( {\frac{n}{2}} \right) + \sum\limits_{l\left| {\left( {\frac{n}{2}} \right)} \right.} {\sigma (l)}  - \sum\limits_{l\left| {\left( {\frac{n}{4}} \right)} \right.} {\sigma (l)} ,\,n \equiv 4\,(\bmod 8) \hfill \\
  \sigma (n) + 3 \cdot \sigma \left( {\frac{n}{2}} \right) + \sum\limits_{l\left| {\left( {\frac{n}{2}} \right)} \right.} {\sigma (l)}  - \sum\limits_{l\left| {\left( {\frac{n}{4}} \right)} \right.} {\sigma (l)}  + 4 \cdot \sum\limits_{l\left| {\left( {\frac{n}{8}} \right)} \right.} {\sigma (l)} \,,\,n \equiv 0\,(\bmod 8) \hfill \\ 
\end{gathered}  \right.$

\item	if $n$ is odd, $c_n = \sigma (n).$

\end {enumerate}
In particular, if $p$ is an odd prime, then $c_p = p + 1.$ 
\end{prop}

\begin{prop}  The number of all normal subgroups of index $n$ in the group $Cc$ is the following:
             \begin{enumerate}	

\item   if $n$ is even, \\
$c_n = \left\{ \begin{gathered}
  \sigma (n) - \sigma \left( {\frac{n}{2}} \right) + \sum\limits_{l\left| {\left( {\frac{n}{2}} \right)} \right.} {\sigma (l)} ,(n \equiv 2 \vee n \equiv 6)\,(\bmod 8) \hfill \\
  \sigma (n) - \sigma \left( {\frac{n}{2}} \right) - \sum\limits_{l\left| {\left( {\frac{n}{4}} \right)} \right.} {\sigma (l)}  + \sum\limits_{l\left| {\left( {\frac{n}{2}} \right)} \right.} {\sigma (l)} ,\,n \equiv 4\,(\bmod 8) \hfill \\
  \sigma (n) - \sigma \left( {\frac{n}{2}} \right) - \sum\limits_{l\left| {\left( {\frac{n}{4}} \right)} \right.} {\sigma (l)}  + \sum\limits_{l\left| {\left( {\frac{n}{2}} \right)} \right.} {\sigma (l)}  + 4 \cdot \sum\limits_{l\left| {\left( {\frac{n}{8}} \right)} \right.} {\sigma (l)} ,\,n \equiv 0\,(\bmod 8) \hfill \\ 
\end{gathered}  \right.$
\item  if $n$ is odd, $c_n = \sigma (n).$
\end {enumerate} 
In particular, if $p$ is an odd prime, then $c_p =p + 1.$ 
\end{prop}

\section{ Proof}

The method for calculating
the number of subgroups of any given index in a group $G$ having an
abelian normal subgroup $T$ of a finite index (see \cite{M}) could be
applied to any polycyclic group or to any poly-(infinite)
cyclic-by-finite group.  In the
sequel, $G$ will denote a space group with the point group isomorphic to
the cyclic group of order 2. In each particular case, we make two standard steps. We firstly
count all subgroups containing $T$ as its subgroup. Then we count all
subgroups contained in $T$.\mbox{} \\

In our setting, $G$ is a finite extension of a free abelian group $T$
of rank 3. Group $T$ is generated by three translations $x$, $y$ and
$z$. A subgroup of finite index in $T$ is free and generated by
elements ${x^a}{y^b}{z^c}$, ${y^d}{z^e}$ and ${z^f}$.
 These exponents are unique to the limits: $a,d,f > 0,0 \le b < d,0 \le c,e < f$.
 Since $P$ is a group isomorphic to the cyclic group of order 2, the index of a subgroup of $T$ in group $G$ is $2adf$.
 We know that the zeta function of $T \cong {\mathbb {Z}^3}$ is $\zeta (s)\zeta (s - 1)\zeta (s - 2)$.
 Therefore, the contribution to the zeta function of a group $G$ coming from this part of the problem
 is ${2^{ - s}}\zeta (s)\zeta (s - 1)\zeta (s - 2)$ for all eight groups.\mbox{} \\ 

Thus, we only need to count subgroups containing $T$. So, let ${H_1}$ be a
subgroup containing  $T$ as its subgroup. Then  ${H_1}$ is generated
by elements $r{x^a}{y^b}{z^c}$, ${x^d}{y^e}{z^f}$, ${y^g}{z^h}$ and
${z^i}$, where $x$, $y$ and $z$ are generators of $T$ and
$rT$ is a generator of the point group $P$, which is cyclic of order 2.
 Furthermore, these exponents are unique to the limits: $d,g,i > 0,0 \le a < d,0 \le b,e < g,0 \le c,f,h < i$.
The index of a subgroup generated by these elements is $dgi$.
Since  $T$ is a normal subgroup in $G$, then ${H_1} \cap T$ is a
normal subgroup in $T$ and \[{{{H_1}} \mathord{\left/
 {\vphantom {{{H_1}} {({H_1} \cap T}}} \right.
 \kern-\nulldelimiterspace} {({H_1} \cap T}}) \simeq {{{H_1}T} \mathord{\left/
 {\vphantom {{{H_1}T} T}} \right. \kern-\nulldelimiterspace} T.}\]
 This means that ${\left( {{x^d}{y^e}{z^f}} \right)^r}$, ${\left( {{y^g}{z^h}} \right)^r}$ and
 ${\left( {{z^i}} \right)^r}$ are elements of ${H_1} \cap T$.
 In this case \[{{{H_1}T} \mathord{\left/ {\vphantom {{{H_1}T} T}} \right. \kern-\nulldelimiterspace} T} \simeq P,\]
 hence ${\left( {r{x^a}{y^b}{z^c}} \right)^2}$ is an element of ${H_1}
\cap T$. The problem of counting subgroups is reduced to solving the
system of equations derived from these conditions. In
effect, we consider the number of possible combinations of values (solutions
of the corresponding system of equations) which the exponents of the
generators of ${H_1}$ may take.\mbox{}
 
If ${H_1}$ is a normal group in $G$,
then its elements have also to satisfy
relations ${\left( {r{x^a}{y^b}{z^c}} \right)^r}$, ${\left(
{r{x^a}{y^b}{z^c}} \right)^x}$, ${\left( {r{x^a}{y^b}{z^c}}
\right)^y}$, ${\left( {r{x^a}{y^b}{z^c}} \right)^z}$ $\in {H_1}$.
Since normality is not a transitive relation, we must also add conditions
which will ensure that a normal subgroup of $T$ is a normal subgroup
of $G$. Let ${H_2} =  < {x^a}{y^b}{z^c},{y^d}{z^e},{z^f} > $ be a
normal group of $T$. Then ${H_2}$ is a normal subgroup of $G$, if
${\left( {{x^a}{y^b}{z^c}} \right)^r}$, ${\left(
{{y^d}{z^e}}\right)^r}$ and $ {\left( {{z^f}} \right)^r}$ are
elements of ${H_2}$.

When writing a space group in an abstract form, we follow the descriptions of
these groups given in \cite{L}. The software packages Mathematica Wolfram and GAP were apt for double checking the 
calculations.  Mathematica was used
to convert our formulas into lists of integers ${a_n}$ or
${c_n}$.

We demonstrate the above technique in detail in the case of group ${C2}$.
\mbox{}

\subsection {\bf Zeta function of group ${C2}$} 

\mbox{}\\

Recall that ${C2} =
\left\langle {\left. {x,y,z,r} \right|\left[ {x,y} \right],\left[
{x,z} \right],\left[ {y,z} \right],{r^2},{x^r} = xy,{y^r} = {y^{ -
1}},{z^r} = {z^{ - 1}}} \right\rangle $.\mbox{} \\

We are counting subgroups of the form ${H_1} = \left\langle
{r{x^a}{y^b}{z^c},{x^d}{y^e}{z^f},{y^g}{z^h},{z^i}} \right\rangle $.
Each of ${\left( {r{x^a}{y^b}{z^c}} \right)^2}$, ${\left( {{x^d}{y^e}{z^f}} \right)^r}$,
   ${\left( {{y^g}{z^h}} \right)^r}$, ${\left( {{z^i}} \right)^r}$ must lie in ${H_1} \cap T$.
Now, ${\left( {r{x^a}{y^b}{z^c}} \right)^2} = {x^{2a}}{y^a}$,
${\left( {{x^d}{y^e}{z^f}} \right)^r} = {x^{ - d}}{y^{ - e}}{z^{ -
f}} \cdot {x^{2d}}{y^d}$, ${\left( {{y^g}{z^h}} \right)^r} = {y^{ -
g}}{z^{ - h}}$, ${\left( {{z^i}} \right)^r} = {z^{ - i}}$. Each of
${x^{ - d}}{y^{ - e}}{z^{ - f}}$, ${y^{ - g}}{z^{ - h}}$,${z^{ -
i}}$ is contained in ${H_1} \cap T = \left\langle
{{x^d}{y^e}{z^f},{y^g}{z^h},{z^i}} \right\rangle $, regardless of
the values of $d$, $e$, $f$, $g$, $h$, $i$. So, ${\left(
{r{x^a}{y^b}{z^c}} \right)^2}$lies in ${H_1} \cap T$ if
${x^{2a}}{y^a}$ lies in ${H_1} \cap T$; ${\left( {{x^d}{y^e}{z^f}}
\right)^r}$ lies in ${H_1} \cap T$ if  ${x^{2d}}{y^d}$ lies in
${H_1} \cap T$. If ${x^{2a}}{y^a}$, ${x^{2d}}{y^d}$ are in ${H_1}
\cap T$ then there exist integers ${\alpha _1},{\alpha _2},{\beta
_1},{\beta _2},{\gamma _1},{\gamma _2}$ such that: \mbox{}
\\${x^{2d}}{y^d} = {\left( {{x^d}{y^e}{z^f}} \right)^{{\alpha
_1}}}{\left( {{y^g}{z^h}} \right)^{{\beta _1}}}{\left( {{z^i}}
\right)^{{\gamma _1}}}$, ${x^{2a}}{y^a} = {\left( {{x^d}{y^e}{z^f}}
\right)^{{\alpha _2}}}{\left( {{y^g}{z^h}} \right)^{{\beta
_2}}}{\left( {{z^i}} \right)^{{\gamma _2}}}$.

We get the following system of equations:

\centerline{ $S_1 = \left\{
\begin{array}{l}
d{\alpha _1} = 2d,e{\alpha _1} + g{\beta _1} = d,f{\alpha _1} + h{\beta _1} + i{\gamma _1} = 0,\\
d{\alpha _2} = 2a,e{\alpha _2} + g{\beta _2} = a,f{\alpha _2} +
h{\beta _2} + i{\gamma _2} = 0,
\end{array} \right\}.$}

The first equation implies ${\alpha _1} = 2$. Consider the equation $2e =
- g{\beta _1} + d$. The left side of the equation is even.  If $g$ is even,
then $d$ has also to be even.  We get:  $\frac{d}{g}
\ge {\beta _1} >  - 2 + \frac{d}{g}$. There are two integers
in the interval $\left[ {\frac{d}{g},\frac{d}{g} - 2} \right)$. To solve
the above system of equations, we shall consider the following cases:

\textbf{Case 1. }  $d$ , $g$, $i$ are odd. In this case, there is
one choice for $a$. Since ${\beta _1}$ has to be odd, there is one
choice for $e$. From $2f + h{\beta _1} + i{\gamma _1} = 0,$ we get $
- \frac{{h{\beta _1}}}{i} \ge {\gamma _1} >  - 2 - \frac{{h{\beta
_1}}}{i}$. Now, $h{\beta _1}$ and $i{\gamma _1}$ have to be odd or even
at the same time. Since ${\beta _1}$ and $i$ are odd, we conclude that if $h$ is
odd then ${\gamma _1}$ is odd and if $h$ is even then ${\gamma _1}$
is even. Hence, there are $i$ choices for $h$ and one choice for
$f$. The zeta function contribution in this case is:
$\sum\limits_{d,g,i \in \mathbb {N}'} {{d^{ - s}}} {g^{ - s}}{i^{ -
s}}gii$, where $\mathbb {N}' = \left\{ {\left. {2k - 1} \right|k \in \mathbb {N}} \right\}$.\mbox{} \\

\textbf{Case 2. and Case 3.} ($d$, $i$ are odd, $g$ is even)  and ( $d$ is odd, $g,i$ are even). These cases are impossible.

\textbf{Case 4.} $i$, $g$ are odd, $d$ is even. There are two
choices for $a$. Since ${\beta _1}$ has to be even, there is one choice
for $e$. Furthermore, ${\gamma _1}$ has to be even, so there is one choice for
$f$. The zeta function contribution in this case is: {$2 \cdot \sum\limits_{d \in 2\mathbb {N},g,i \in \mathbb {N}'} {{d^{ - s}}}
{g^{ -s}}{i^{ - s}}gii$, where $\mathbb {N}' = \left\{ {\left. {2k - 1} \right|k \in \mathbb {N}} \right\}$.\mbox{} \\

\textbf{Case 5.} $g$, $d$ are odd, $i$ is even. There is one
choice for $a$; ${\beta _1}$ has to be odd, so there is one choice
for $e$. Since $h{\beta _1}$ and $i{\gamma _1}$ have to be odd or
even at the same time, we see that $h$ has to be even. There are two
choices for ${\gamma _1}$, hence there are two choices for $f$. The zeta
function contribution in this case is: $ \sum\limits_{i \in 2\mathbb {N},g,d \in \mathbb {N}'} {{d^{ - s}}}
{g^{ -s}}{i^{ - s}}gi \cdot \frac{i}{2} \cdot 2$. \mbox{} \\

\textbf{Case 6.} $g$ is odd, $i$, $d$ are even. There are two choices for $a$; ${\beta _1}$ has to be
even, so there is one choice for $e$. If $a = 0$, then there are two
choices for $f$; if $a = \frac{d}{2}$ then ${\gamma _1}$ is even, so
 there is one choice for $f$ in this case.
The zeta function contribution reads: $3\sum\limits_{g \in \mathbb {N}',d,i \in 2\mathbb {N}} {{d^{ - s}}} {g^{ - s}}{i^{ - s}}gii$. \mbox{} \\

\textbf{Case 7.} $g$, $d$ are even, $i$ is odd. There are two choices for $a$. If $a = 0$ and ${\beta _1}$
is  even, then ${\gamma _1}$ is even and there are $i$ choices for
$h$. Hence, for $a = 0$ and ${\beta _1}$ is even, there are one
choice for $e$ and one choice for $f$. If $a = 0$ and ${\beta _1}$
is  odd, then ${\gamma _1}$ and  $h$ are both even or both odd.
If $a = \frac{d}{2}$ then ${\beta _1}, {\gamma _1}$
are even and there
is one choice for $f$.
The zeta
function contribution in this case is: $3\sum\limits_{i \in \mathbb {N}',d,g \in 2\mathbb {N}} {{d^{ - s}}} {g^{ - s}}{i^{ - s}}gii$.\mbox{} \\

\textbf{Case 8.} $g,d$ $i$ are even. There are two choices for $a$. If $a = 0$ and ${\beta _1}$
is  even, then there are two choice for  ${\gamma _1}$ and there are
$i$ choices for $h$. If $a = 0$ and ${\beta _1}$ is odd, then there
are two choices for ${\gamma _1}$ and  $h$ is even ($\frac{i}{2}$
choices for $h$). If $a = \frac{d}{2}$ then ${\beta _1}$, ${\gamma
_1}$ are even and there is one choice for $f$.
The zeta
function contribution in this case is:
\mbox{}\\$\sum\limits_{d,i,g \in 2\mathbb {N}} {{d^{ - s}}} {g^{ - s}}{i^{ - s}} \cdot g \cdot i \cdot i \cdot 2 + \sum\limits_{d,i,g \in 2\mathbb {N}} {{d^{ - s}}} {g^{ - s}}{i^{ - s}} \cdot g \cdot i \cdot \frac{i}{2} \cdot 2 + \sum\limits_{d,i,g \in 2\mathbb {N}} {{d^{ - s}}} {g^{ - s}}{i^{ - s}} \cdot g \cdot i \cdot i = 4\sum\limits_{d,i,g \in 2\mathbb {N}} {{d^{ - s}}} {g^{ - s}}{i^{ - s}} \cdot g \cdot i \cdot i$. \mbox{} \\

Finally, we obtain the zeta function for group  ${C2}$:\mbox{} \\

${\zeta _{{C2}}}(s) = \sum\limits_{d,g,i \in \mathbb {N}'} {{d^{ - s}}} {g^{
- s}}{i^{ - s}}gii + 2 \cdot \sum\limits_{d \in 2\mathbb {N},g,i \in \mathbb {N}'} {{d^{
- s}}} {g^{ - s}}{i^{ - s}}gii +\\+ \sum\limits_{d,g \in \mathbb {N}',i \in 2\mathbb {N}}
{{d^{ - s}}} {g^{ - s}}{i^{^{ - s}}}gii +3\sum\limits_{g \in \mathbb {N}',i,d \in 2\mathbb {N}} {{d^{ - s}}} {g^{ - s}}{i^{ - s}}gii + 3\sum\limits_{i \in \mathbb {N}',d,g \in 2\mathbb {N}} {{d^{ - s}}} {g^{ - s}}{i^{ - s}}gii +\\+
4\sum\limits_{d,i,g \in 2\mathbb {N}} {{d^{ - s}}} {g^{ - s}}{i^{ - s}}gii +
{2^{ - s}}\zeta (s)\zeta (s - 1)\zeta (s - 2) = (1 + {2^{ - 2s +
3}})\zeta (s)\zeta_1 (s)\zeta_2 (s)$.

\subsection{\bf Normal zeta function of group ${C2}$}
\mbox{}\\

We use the set of
constraints which we obtained in the previous Subsection. By counting the number of subgroups ${H_1} = \left\langle
{r{x^a}{y^b}{z^c},{x^d}{y^e}{z^f},{y^g}{z^h},{z^i}} \right\rangle $
of ${C_2}$}, we deduced the system ${S_1}$.

Based on the conditions of normality, we get another set
of constraints: \\
 \\\centerline{${S_2} = \left\{
\begin{array}{l}
d{\alpha _3} = 0,e{\alpha _3} + g{\beta _3} =  - 1,f{\alpha _3} + h{\beta _3} + i{\gamma _3} = 0\\
d{\alpha _4} = 0,e{\alpha _4} + g{\beta _4} = 2,f{\alpha _4} + h{\beta _4} + i{\gamma _4} = 0,\\
d{\alpha _5} = 0,e{\alpha _5} + g{\beta _5} = 0,f{\alpha _5} +
h{\beta _5} + i{\gamma _5} = 2
\end{array} \right\}.$}\\

The equations $d{\alpha _3} = 0,e{\alpha _3} + g{\beta
_3} =  - 1,f{\alpha _3} + h{\beta _3} + i{\gamma _3} = 0$ imply that $g = 1,{\beta _3} = -1$, so $b = e = 0 = h$,
while the equations  $d{\alpha _5} = 0,e{\alpha _5} + g{\beta _5} =
0,f{\alpha _5} + h{\beta _5} + i{\gamma _5} = 2$ imply that $i = 1$
or $i = 2$. Observing four cases depending on whether $d$ is even
or odd and depending on values of  $i$, we get
 \mbox{} \\ $\sum\limits_{d
\in \mathbb {N}'} {{d^{ - s}}}  + 2\sum\limits_{d \in 2\mathbb {N}} {{d^{ - s}}}  + 4
\cdot {2^{ - s}}\sum\limits_{d \in \mathbb {N}} {{d^{ - s}}}  +
6\sum\limits_{d \in 2\mathbb {N}} {{d^{ - s}}} {2^{ - s}} = (2 \cdot {2^{ -
2s}} + 5 \cdot {2^{ - s}} + 1)\zeta (s)$. \\

Now, we count normal subgroups of $T$. Any such subgroup takes
the form ${H_2}$$ = \left\langle {{x^a}{y^b}{z^c},{y^d}{z^e},{z^f}}
\right\rangle $. We assume $0 < a,0 \le b < d,0 \le c,e <
f$. Let us check the conditions of normality in ${C2}$. These require
that ${\left( {{x^a}{y^b}{z^c}} \right)^r}$, ${\left( {{y^d}{z^e}}
\right)^r}$, ${\left( {{z^f}} \right)^r}$ are in ${H_2}$. After
some calculations, we get the next set of
constraints:
 \\\centerline{${S_3} = \left\{ {a{\alpha _1} = 2a,b{\alpha _1} +
d{\beta _1} = a,c{\alpha _1} + e{\beta _1} + f{\gamma _1} = 0}
\right\}.$}\\

This system of equations reduces to three equations:
$d{\beta _1} = a - 2b,e{\beta _1} + f{\gamma _1} = 0, - 2e +
f{\gamma _1} = 0$. The case $a$ odd and $d$ even is
impossible. Thus, we have six cases depending on whether $a$ or $d$ or
$f$ are odd or even. The respective contributions are:\\

$\sum\limits_{a,d,f \in \mathbb {N}'} {{a^{ - s}}} {d^{ - s}}{f^{ -
s}}f$ $ + \sum\limits_{a \in 2\mathbb {N},d,f \in \mathbb {N}'} {{a^{ - s}}} {d^{ -
s}}{f^{ - s}}f$$ + 2\sum\limits_{a,f \in 2\mathbb {N},d \in \mathbb {N}'} {{a^{ - s}}}
{d^{ - s}}{f^{ - s}}f$$ +\\+2\sum\limits_{a,d \in 2\mathbb {N},f \in \mathbb {N}'}
{{a^{ - s}}} {d^{ - s}}{f^{ - s}}f$+$\sum\limits_{a,d \in \mathbb {N}',f \in 2\mathbb {N}} {{a^{ - s}}} {d^{ - s}}{f^{ - s}}f + 3\sum\limits_{a,d,f \in 2\mathbb {N}}
{{a^{ - s}}} {d^{ - s}}{f^{ - s}}f$$ \\ = (1 - {2^{ - s}} + 4 \cdot
{2^{ - 2s}})\zeta (s)\zeta (s)\zeta (s - 1)$.\\

Combining the contributions coming from ${H_1}$ and ${H_2}$, we
get the normal zeta function of ${C2}$: \mbox{} \\ $\zeta _{{C2}}^
\triangleleft (s) = (2 \cdot {2^{ - 2s}} + 5 \cdot {2^{ - s}} +
1)\zeta (s) + {2^{ - s}} \cdot (1 - {2^{ - s}} + 4 \cdot {2^{ -
2s}})\zeta (s)\zeta (s)\zeta_1 (s)$.

\section {\bf Remarks}

As far as subsection 5.1  is concerned, using a similar argumentation one can obtain 
more general assertions of the following form.

\begin{prop}  Let $n \ge 2 $ be a fixed integer and ${G_n}$ be the
group defined by \mbox{} \\ ${G_n} = \left\langle {\left.
{r,{x_1},{x_2},...,{x_n}} \right|{r^2},\left[ {{x_j},{x_k}}
\right](\forall 1 \le j,k \le n),{x_i}^r = {x_i}^{ - 1}(\forall 1
\le i \le n)} \right\rangle $.
  The zeta function of group  ${G_n}$ is given by\\
\centerline{${\zeta _{{G_n}}}(s) = \zeta (s - 1)\zeta (s - 2) \cdots
\zeta (s - n) + {2^{ - s}}\zeta (s)\zeta (s - 1)\zeta (s - 2) \cdots
\zeta (s - n +1).$}
\end{prop}

\begin{prop} Let   $n \ge 2 $ be an integer and ${G_n}$ be the group
defined by\\ ${G_n}=\left\langle {\left. {r,{x_1},{x_2},...,{x_n}} \right|{r^2},\left[ {{x_j},{x_k}} \right]
  (\forall 1 \le j,k \le n),{x_1}^r = {x_1}{x_2},{x_i}^r={x_i}^{ - 1}(\forall 1< i \le n)} \right\rangle $.
The zeta function of group ${G_n}$ is
 ${\zeta _{{G_n}}}(s) = (1 +{2^{ - 2s + n}})\zeta (s)\zeta (s - 1)\zeta (s - 2) \cdots \zeta (s
- n + 1)$. The zeta function of group ${G_n}$  has an Euler product.
\end {prop}

\end{document}